\documentclass[11pt]{article}

\usepackage{latexsym}
\usepackage{amsfonts}
\usepackage{amsmath}
\newtheorem{theorem}{Theorem}[section]

\newtheorem{proposition}{Proposition}[section]

\begin{document}

\title{Azuma-Hoeffding bounds for a class of urn models}

\author{Amites Dasgupta\thanks{%Postal Address:
Stat-Math Unit, Indian Statistical Institute, 203 B. T. Road,
Kolkata 700108, INDIA. Email Address: amites@isical.ac.in} \\
Indian Statistical Institute}
%\date{}

\maketitle

\begin{abstract}

We obtain Azuma-Hoeffding bounds, which are exponentially decreasing,
for the probabilities of being away from the limit
for a class of urn models. The method consists of relating the variables to certain
linear combinations using eigenvectors of the replacement matrix, thus bringing in
appropriate martingales. Some cases of repeated eigenvalues are also considered using
Jordan vectors.

\end{abstract}


{\bf 2000 Subject classification:} Primary: 60F10, 60G42, Secondary: 60C05, 60E15. 

\section{Introduction}

Consider a two color urn model, with colors white and black, evolving as follows: at time
$0$ the color count is $(W_0, B_0) = \boldsymbol{C}_0, 0 \leq W_0, B_0 \leq 1, W_0 + B_0 = 1$.
There is a $2 \times 2$ irreducible and stochastic replacement matrix $R$ which 
drives the color count evolution as follows: given the color composition $(W_n, B_n) = \boldsymbol{C}_n$
at time $n$, we choose a color at random with probabilities proportional to $W_n/(n+1), B_n/(n+1)$,
respectively. If white is chosen we add the first row of $R$, or if black is chosen we add the 
second row of $R$, to $\boldsymbol{C}_n$ to get $\boldsymbol{C}_{n+1}$, the color count vector at time $(n+1)$. Denoting by
$\boldsymbol{\chi}_{n+1}$ a row vector which takes the value $(1, 0)$ if white is chosen, or $(0, 1)$
if black is chosen as above, we clearly have at the $(n+1)$th stage $\boldsymbol{C}_{n+1} =
\boldsymbol{C}_{n} + \boldsymbol{\chi}_{n+1} R$. Since each row sum of $R$ is one, the total color count
in each step increases by 1, making the urn model balanced.Various eigenvalues and eigenvectors of $R$ 
can be used to understand the limiting behavior of various linear combinations of $W_n$ and $B_n$. The above model easily
generalizes to urn models with more than two colors and we state some of the limit theorems (for these more
general models with $d$ colors) necessary for our notation and results. 

The strong law states that $\frac{\boldsymbol{C}_n}{n+1} \rightarrow \boldsymbol{\pi}$ almost surely, where $\boldsymbol{\pi}$
is the left eigenvector to the Perron-Frobenius eigenvalue 1 of $R$ (see Laruelle and Pages \cite{laruellepages}).
This $\boldsymbol{\pi}$ is also a probability vector with all components positive. Next, if
we multiply the equation $\boldsymbol{C}_{n+1} = \boldsymbol{C}_{n} + \boldsymbol{\chi}_{n+1} R$ by a right 
eigenvector $\boldsymbol{\xi}$ to a non Perron-Frobenius eigenvalue $\lambda$ (which necessarily satisfies $|\lambda| < 1$),
then $R \boldsymbol{\xi} = \lambda \boldsymbol{\xi}$ gives a reduction which has been used in central limit theorems.
The scaling in these theorems depend on whether the real part of $\lambda$ is $1/2$ or not. See Smythe \cite{smythe}, Basak and Dasgupta 
\cite{bd}, and the comprehensive paper Janson \cite{janson}.   

The above two results are clearly the analogues (for urn models) of the strong law of large numbers and the central
limit theorem for the sum $S_n$ of i.i.d. random variables $X_i$ with mean $\mu$ and variance $\sigma^2$. For sums of i.i.d.
random variables, another class of 
very important reasults measure the probability of deviation $P(|\frac{S_n}{n} - \mu| > \epsilon) = P(|S_n - nu| > \epsilon n)$.
If these probabilities decrease fast (say like $e^{-const. n}$) then even for $n$ not too large, $S_n/n$ is close to $\mu$
with high probability. In fact, under the assumption of finite exponential moments, in the theory of large deviation 
the limit $\lim \frac{1}{n} \log P(|\frac{S_n}{n} - \mu| > \epsilon)$ is obtained. However in the general dependent set up
such a limit is not easy to obtain (see however Grama and Haeusler \cite{grama}, for an extension), and for the required random variables 
$Y_n/n$ with limit $\nu$ say, one sided
bounds like $P(|Y_n - n \nu| > n \epsilon) \leq e^{-const. f(n)}$ where $f(n)$ is an increasing function, serve the same
useful purpose. One such bound is the Azuma inequality, also called Azuma-Hoeffding bound, see Ross \cite{ross}, which
assumes a martingale structure.

In the case of our urn models clearly the last problem reduces to finding bounds for $P(|W_n - (n + 1) \pi_1| > (n+1) t)$. 
Because of the dependence here, large deviation considerations have been approached in different ways in the literature. 
For an approach through generating functions see Flajolet, Dumas, Puyhaubert \cite{flajolet} and Morcrette \cite{morcrette}.
Franchini \cite{franchini} gives a functional form for a class of urn models, Bandyopadhyay and Thacker \cite{thacker} 
consider deviations from the expected configuration over a class of probabilities which keep the expected configuration
fixed and get rate $\log n$. Even in a simple set up, it is not easy to understand the limit (or limsup and liminf) of 
$- \frac{1}{n} \log P(\boldsymbol{C}_n/(n+1) \in B)$, where $\boldsymbol{C}_n$ is the urn composition at time $n$, and
$B$ is a set in $\mathbb{R}^d$ where $d$ is the number of colors. We thus focus on exponential upper bounds depending on the eigenvalues for the 
class of urn models described in the first paragraph.  
 
We now describe the organization of the article. In section 2 we consider a linear combination using the eigenvector
$\boldsymbol{\xi}$ to the non Perron-Frobenius eigenvalue $\lambda$ and show how deviations of $W_n$ can be related
to deviations of a corresponding martingale. In section 3 the same is done for some multicolor urn models, however
now there may be more than one distinct eigenvalues, or a repeated eigenvalue. We consider only the case when these eigenvalues
are real and stress the systematic use of eigenvectors and Jordan vectors.

\section{Two Colors}

Let us now consider a  two color urn model evolving following a $2 \times 2$ irreducible stochastic replacement
matrix $R$. The left eigenvector $\boldsymbol{\pi} = (\pi_1, \pi_2)$ corresponds to the Perron-Frobenius eigenvalue
(will be called principal eigenvalue from now on) 
1, the right eigenvector $\boldsymbol{\xi} = (\xi_1, \xi_2)^T$ corresponds to the other nonprincipal eigenvalue
$\lambda \in (-1, 1)$. Let $\boldsymbol{C}_n = (W_n, B_n)$ denote the composition at time $n$, it
being assumed that $W_0 + B_0 = 1$, so that $W_n + B_n = n + 1$. It is known from the strong law that
$\boldsymbol{C}_n/(n+1) \rightarrow \boldsymbol{\pi}$ almost surely. In addition $1. \boldsymbol{\pi} . \boldsymbol{\xi}
= \boldsymbol{\pi} R \boldsymbol{\xi} = \lambda \boldsymbol{\pi}. \boldsymbol{\xi}$ and $\lambda \in (-1, 1)$ 
implies $\pi_1 \xi_1 + \pi_2 \xi_2 = 0$.

Suppose we can derive an inequality for $P(\boldsymbol{C}_n . \boldsymbol{\xi} > (n+1)t)$ where $t > 0$. Notice that
$W_n \xi_1 +  (n+1 - W_n)\xi_2 > (n+1)t$ iff $W_n (\xi_1 - \xi_2) > - (n + 1) \xi_2 + (n+1)t$. But $\pi_1 + \pi_2 = 1,
\pi_1 \xi_1 + \pi_2 \xi_2 = 0$, makes the preceding inequality $W_n > (n+1) \pi_1 + (n+1)(t/(\xi_1 - \xi_2))$ 
assuming $\xi_1 > 0, \xi_2 < 0$. Thus probability inequalities for $W_n$ can also be obtained if we obtain
probability inequalities for linear combinations corresponding to an eigenvector of the nonprincipal eigenvalue.

In view of this we represent $\boldsymbol{C}_n . \boldsymbol{\xi}$ as a sum of martingale differences and apply
Azuma-Hoeffding inequality to get exponential bounds. In many practical applications such bounds suffice, since large deviation principle (LDP) 
under dependent set up is not easy. The exponential bounds are well known in the urn model literature for particular examples, 
also see the recent work of Kuba and Sulzbach \cite{kuba} who also use martingale inequalities for Polya type urns, where we note that
our irreducible $R$ is different from the $R$ of Polya's original urn models. Our approach is through a combination of linear algebra 
and martingale theory which tries to cover a family of urn models depending on the value of $\lambda$. Although
the precise forms can be recovered from the proof that follows, the statement of next proposition stresses
the increasing nature of $f(n)$. 

We first state the form of Azuma-Hoeffding inequality that we use. Suppose $(M_i, {\mathcal F}_i)_{0 \leq i \leq n}$ 
is a martingale such that the increments $\Delta M_i = M_{i+1} - M_i$ are bounded as follows, $- \alpha_i \leq \Delta M_i \leq \beta_i,
0 \leq i \leq n-1$. The Azuma inequality is the following exponential inequality for the probability of large deviations
$P(M_n - M_0 > tn) \leq \exp \{ - 2n^2t^2/(\sum (\alpha_i + \beta_i)^2) \}, t>0$, with similar statements for $t < 0$.
Our martingale differences depend on both $i$ and $n$, and handling this is the main part of the work.

\begin{theorem}\label{prop1}
For $\lambda \in (-1, 1)$ one can find increasing functions $f(n)$ so that for positive $t$,
\[ P(\boldsymbol{C}_{n+1} . \boldsymbol{\xi} - \Big{(} \prod_0^n (1 + \frac{\lambda}{j+1}) \Big{)} \boldsymbol{C}_0.\boldsymbol{\xi} > 
(n+1)t) \leq  e^{- c. t^2 f(n)},\] where $c$ is a constant greater than zero. 
\end{theorem}

{\bf Proof:} When $\lambda = 0$ the martingale doesn't move. Since $R \boldsymbol{\xi} = \boldsymbol{0}$ we have 
$\boldsymbol{C}_{n+1} \boldsymbol{\xi} = \boldsymbol{C}_n
\boldsymbol{\xi}$, that is this linear combination equals $\boldsymbol{C}_0 \boldsymbol{\xi}$ and the exponential 
inequality is satisfied for any $c > 0$.

For $\lambda \neq 0$, consider the equations
\begin{eqnarray}\label{evolution}
\boldsymbol{C}_{n+1} . \boldsymbol{\xi} &=& \boldsymbol{C}_n .\boldsymbol{\xi} + \lambda \boldsymbol{\chi}_{n+1}. \boldsymbol{\xi} \nonumber \\
&=& \boldsymbol{C}_n .\boldsymbol{\xi} + \lambda \frac{\boldsymbol{C}_n . \boldsymbol{\xi}}{n+1} + \lambda (
\boldsymbol{\chi}_{n+1}. \boldsymbol{\xi} - \frac{\boldsymbol{C}_n . \boldsymbol{\xi}}{n+1}) \nonumber \\
&=& (1 + \frac{\lambda}{n+1}) \boldsymbol{C}_n . \boldsymbol{\xi} + \lambda (
\boldsymbol{\chi}_{n+1}. \boldsymbol{\xi} - \frac{\boldsymbol{C}_n . \boldsymbol{\xi}}{n+1}),
\end{eqnarray}
the martingale differences coming from $E(\boldsymbol{\chi}_{n+1}|\mathcal{F}_n) = \boldsymbol{C}_n/(n+1)$ where $\mathcal{F}_n$
is the urn composition $\sigma$-field upto time $n$. This gives us 
the iteration 
\begin{equation}\label{martingale}
\boldsymbol{C}_{n+1} . \boldsymbol{\xi} = \Big{(} \prod_0^n (1 + \frac{\lambda}{j+1}) \Big{)} \boldsymbol{C}_0.\boldsymbol{\xi} +
\sum_{j = 0}^n \Big{\{} \prod_{k = j+1}^n (1 + \frac{\lambda}{k+1}) \Big{\}} 
\lambda (
\boldsymbol{\chi}_{j+1}. \boldsymbol{\xi} - \frac{\boldsymbol{C}_j , \boldsymbol{\xi}}{j+1}),
\end{equation}
where we use the notational convention that for $j = n$, the product \\ $\prod_{k = j+1}^n (1 + \frac{\lambda}{k+1})$ will be understood
to be 1.

 Since $\boldsymbol{\chi}_{n+1}$ and 
$\boldsymbol{C}_n/(n+1)$ have bounded components, $\lambda (
\boldsymbol{\chi}_{j+1}. \boldsymbol{\xi} - \frac{\boldsymbol{C}_j . \boldsymbol{\xi}}{j+1})$ is bounded by the same constant for all $j$,
and to understand the exponent in the Azuma inequality we need to estimate $$(n + 1)^2/ \sum_{j = 0}^n \Big{\{} \prod_{k = j+1}^n (1 +
\frac{\lambda}{k+1}) \Big{\}}^2 .$$
We further notice that if 
\begin{equation}\label{dn} 
\sum_{j = 0}^n \Big{\{} \prod_{k = j+1}^n (1 +
\frac{\lambda}{k+1}) \Big{\}}^2 \leq D_{n}(\lambda),
\end{equation} 
then $e^{- 2(n + 1)^2t^2/ \sum_{j = 1}^n \Big{\{} \prod_{k = j+1}^n (1 +
\frac{\lambda}{k+1}) \Big{\}}^2} \leq e^{- 2(n+1)^2 t^2/D_{n}(\lambda)}$. Hence we look for successive upper bounds for the sum.  
First using Euler's relation $\prod_0^n (1 + \frac{\lambda}{j+1}) \sim \frac{n^{\lambda}}{\Gamma(\lambda + 1)}$, the term for 
$j=0$ is $\sim n^{2\lambda}$, and 
the products $\Big{\{} \prod_{k = j+1}^n (1 +
\frac{\lambda}{k+1}) \Big{\}}^2, j = 1, \ldots, n$, are bounded above by $const. (n/j)^{2\lambda}$, the constant being uniform over $n$ and $j$. Next to bound the sum 
$\sum_1^n (1/j)^{2\lambda}$ from above we note that the functions $g(x) = 1/x^{2\lambda}$ are increasing for $\lambda < 0$
and decreasing for $\lambda >0$. Using Euler's comparison between sums and integrals for such functions we
have for increasing $g$, $g(1) + \cdots + g(n) \leq \int_1^{n+1} g(x) dx$ and for decreasing $g$, $g(1) + \cdots + g(n)
\leq g(1) + \int_1^n g(x) dx$.

In the asymptotics $\int_1^n g(x) dx$ and $\int_1^{n+1} g(x) dx$ do not make any difference, thus we simply
consider the behavior of $n^{2\lambda} \int_1^n (1/x^{2\lambda}) dx$. The behaviors in the three cases are 
(a) $n - n^{2\lambda} \sim n$ for $\lambda < 1/2 (\lambda \neq 0)$, (b) $n \log n$ for $\lambda = 1/2$, and (c) $n^{2\lambda} - n \sim n^{2\lambda}$
for $1/2 < \lambda < 1$. Also for $\lambda > 0$ the behaviors of $n^{2\lambda} g(1)$ are like (a) $o(n), \lambda < 1/2$,
(b) $n, \lambda = 1/2$, (c) $n^{2\lambda}, \lambda > 1/2$. Thus for the upper bound of the sum  we have
$D_{n}(\lambda) \sim$ (a) $n, \lambda < 0$, (b) $n + o(n), 0 < \lambda < 1/2$, (c) $n + n \log n, \lambda = 1/2$,
(d) $n^{2\lambda} + n^{2\lambda}, 1/2 < \lambda < 1$, including the term for $j = 0$ which is $\sim n^{2\lambda}$.
Since the bound from the Azuma inequality 
is  bounded by $e^{-2(n+1)^2 t^2/D_{n}(\lambda)}$, we get the increasing functions $f(n)$ in the statement 
of the proposition. \hfill $\Box$

{\bf Remark:} The $0$-th term $A_n(\lambda) = \Big{(} \prod_0^n (1 + \frac{\lambda}{j+1}) \Big{)} \boldsymbol{C}_0.\boldsymbol{\xi}$
is $O(n^\lambda), - 1 < \lambda < 1$, and $\boldsymbol{C}_0$ is constant. For large $n$, $n$ dominates $n^\lambda$ and since
the limit of $\boldsymbol{C}_{n+1} \boldsymbol{\xi}/(n+1)$ is zero almost surely, for large deviation purposes the regions
$\boldsymbol{C}_{n+1} \boldsymbol{\xi}/(n+1) - A_n/(n+1)  > t_1$ and $ \boldsymbol{C}_{n+1} \boldsymbol{\xi}/(n+1) > 
t_2$ can be compared and a slightly different region obtained for the latter depending on $n$.

\section{Some multicolor cases with real distinct, or repeated eigenvalues}

With the notation as before, consider a three color urn for which $R$ is a $3 \times 3$ irreducible stochastic matrix, $\boldsymbol{\pi}$ is the left
eigenvector to the principal eigenvalue 1. Now consider the right eigenvector $\boldsymbol{1}$ to eigenvalue $1$ consisting of all 1's, and
suppose there are two other linearly independent right eigenvectors $\boldsymbol{\xi}_2$ and $\boldsymbol{\xi}_3$ corresponding
to real distinct nonprincipal eigenvalues $\lambda_2, \lambda_3$ respectively. In order to get the count of the first color say,
we need to multiply $\boldsymbol{C}_n$ by $(1, 0, 0)^T$. Now by linear independence $(1, 0, 0)^T = \alpha_1 \boldsymbol{1}
+ \alpha_2 \boldsymbol{\xi}_2 + \alpha_3 \boldsymbol{\xi}_3$ for some $\alpha_i, i = 1, 2, 3$. Since total color count at time $n$ is $(n+1)$ and
$\boldsymbol{C}_n.\boldsymbol{\xi}_i/(n+1) \rightarrow \boldsymbol{\pi} \boldsymbol{\xi}_i = 0$, almost surely, $\alpha_1 = \pi_1$. Thus $C_{1n} - 
\pi_1 (n+1) = \alpha_2 \boldsymbol{C}_n. \boldsymbol{\xi}_2 + \alpha_3 \boldsymbol{C}_n. \boldsymbol{\xi}_3$ and similarly for other colors.
In view of this $P(C_{1n} - \pi_1 (n+1) > t(n+1)) = P(\alpha_2 \boldsymbol{C}_n. \boldsymbol{\xi}_2 
+ \alpha_3 \boldsymbol{C}_n . \boldsymbol{\xi}_3 > t(n+1))$. In this case the two martingale differences corresponding to the 
two eigenvectors can be added and we can use the calculations of the previous proof as follows (of the two distinct ones, one zero eigenvalue  
contributes a constant $0$-th term only):

\begin{theorem}\label{prop2}
For nonzero $\lambda_2, \lambda_3 \in (-1, 1)$ one can find increasing functions $f(n)$ so that for positive $t$,
\begin{eqnarray*} && P(\alpha_2 \boldsymbol{C}_{n+1} . \boldsymbol{\xi_2} + \alpha_3 \boldsymbol{C}_{n+1} . \boldsymbol{\xi_3} 
 - \alpha _2 \Big{(} \prod_0^n (1 + 
\frac{\lambda_2}{j+1}) \Big{)} \boldsymbol{C}_0.\boldsymbol{\xi_2} \\
&& - \alpha _3 \Big{(} \prod_0^n (1 + 
\frac{\lambda_3}{j+1}) \Big{)} \boldsymbol{C}_0.\boldsymbol{\xi_3} > 
(n+1)t) \\ &\leq & e^{- c. t^2 f(n)},
\end{eqnarray*}
where $c$ is a constant greater than zero. 
\end{theorem}

{\bf Proof:} The combined martingale differences here give rise to the martingale difference 
\begin{eqnarray*}
&& \alpha_2 \Big{\{} \prod_{k = j+1}^n (1 + \frac{\lambda_2}{k+1})
\Big{\}}  \lambda_2 (
\boldsymbol{\chi}_{j+1}. \boldsymbol{\xi}_2 - \frac{\boldsymbol{C}_j . \boldsymbol{\xi}_2}{j+1}) \\
&+& \alpha_3 \Big{\{} \prod_{k = j+1}^n 
(1 + \frac{\lambda_3}{k+1}) \Big{\}} 
\lambda_3 (
\boldsymbol{\chi}_{j+1}. \boldsymbol{\xi_3} - \frac{\boldsymbol{C}_j . \boldsymbol{\xi_3}}{j+1}).
\end{eqnarray*}
As before these martingale differences are bounded by constant times $$\prod_{k = j+1}^n (1 + \frac{\lambda_2}{k+1}) + \prod_{k = j+1}^n (1 + \frac{\lambda_3}{k+1}),$$
which are bounded by constant times $\frac{n^{\lambda_2}}{j^{\lambda_2}} + \frac{n^{\lambda_3}}{j^{\lambda_3}}$. For the squares of these bounds 
we can use $(a + b)^2 \leq 2 (a^2 + b^2)$ and repeat the previous proof starting analogously from
inequality (\ref{dn}) with $D_{n}(\lambda_2) + D_{n}(\lambda_3)$. \hfill $\Box$

Now suppose that $R$ has one repeated real nonprincipal eigenvalue $\lambda$. Using $R P = P \begin{pmatrix}
\lambda & 1 \\ 0& \lambda \end{pmatrix}$ where the columns of $P$ are $\boldsymbol{\xi}_2, \boldsymbol{\xi}_3$
respectively, we get $R \boldsymbol{\xi}_3 = \boldsymbol{\xi}_2 + \lambda
\boldsymbol{\xi}_3$. Like the eigenvector $\boldsymbol{\xi}_2$, for the Jordan vector $\boldsymbol{\xi}_3$  one again has 
$\boldsymbol{C}_n \boldsymbol{\xi}_3/(n+1) \rightarrow 0$. To see this note that with ${\mathcal F}_n$, the urn
composition $\sigma$-field upto time $n$, $E(\boldsymbol{C}_{n+1} \boldsymbol{\xi}_3|{\mathcal F}_n) = (1 + \frac{\lambda}{n+1}) \boldsymbol{C}_n 
\boldsymbol{\xi}_3 + \frac{\boldsymbol{C}_n}{n+1} \boldsymbol{\xi}_2$. In this case following Dasgupta and Maulik \cite{dasguptamaulik} one can 
consider the martingale \[ M_n = \frac{\boldsymbol{C}_n \boldsymbol{\xi}_3}{\Pi_n(\lambda)} - \sum_{j = 0}^{n-1}
\frac{1}{(j+1) \Pi_{j+1}(\lambda)} \boldsymbol{C}_j \boldsymbol{\xi}_2,  \] with $\Pi_n(\lambda) = \prod_{j=0}^{n-1}
(1 + \frac{\lambda}{j+1})$. The martingale differences of the above $M_n$ have variances like $1/n^{2 \lambda}$
and Lemma 2.1 of Dasgupta and Maulik \cite{dasguptamaulik} with $a_n = n^{1-\lambda}$ shows $\boldsymbol{C}_n \boldsymbol{\xi}_3/(n+1) \rightarrow 0$
using $\boldsymbol{C}_j \boldsymbol{\xi}_2/(j + 1) \rightarrow 0$.
Martingales with other Jordan vectors can be handled in a similar manner successively.
Hence as before, for $\boldsymbol{C}_{1n} - \pi_1(n+1)$, it is enough to consider a linear combination using 
the eigenvector $\boldsymbol{\xi}_2$ and the Jordan vector $\boldsymbol{\xi}_3$. From  the statements of
the previous theorems and the proofs we see that
identifying the $0$-th term and the sum of the martingale differences for each linear combination is necessary. We first do this for
the case of the repeated eigenvalue zero to bring out the different nature of $\boldsymbol{C}_{n+1} \boldsymbol{\xi}_3$ 
corresponding to the Jordan vector.

\begin{proposition}
For the repeated eigenvalue $\lambda = 0$,  $\boldsymbol{C}_{n+1} \boldsymbol{\xi}_2 = \boldsymbol{C}_0 \boldsymbol{\xi}_2$ 
and $\boldsymbol{C}_{n+1} \boldsymbol{\xi}_3$ has a $0$-th term which is $o(n)$ and martingale differences which are
bounded.
\end{proposition}

{\bf Proof:} From the proof of Theorem \ref{prop1}, $\boldsymbol{C}_{n+1} \boldsymbol{\xi}_2$ is constant ($= \boldsymbol{C}_0 \boldsymbol{\xi}_2$)
since $\boldsymbol{\xi}_2$ is the eigenvector to eigenvalue $0$. Next, since $R \boldsymbol{\xi}_3 = \boldsymbol{\xi}_2 + 0 . \boldsymbol{\xi}_3$, we have
$\boldsymbol{C}_{n+1} \boldsymbol{\xi}_3 = \boldsymbol{C}_n \boldsymbol{\xi}_3 + \boldsymbol{\chi}_{n+1} \boldsymbol{\xi}_2 = \boldsymbol{C}_n \boldsymbol{\xi}_3 
+ \frac{1}{n+1} \boldsymbol{C}_n \boldsymbol{\xi}_2 + ( \boldsymbol{\chi}_{n+1} - \frac{\boldsymbol{C}_n}{n+1})\boldsymbol{\xi}_2$.
 
Iterating this equation $\boldsymbol{C}_{n+1} \boldsymbol{\xi}_2 = \boldsymbol{C}_0 \boldsymbol{\xi}_3 + 
\boldsymbol{C}_0  \boldsymbol{\xi}_2 \sum_{j=0}^n \frac{1}{j+1} +
\sum_{j = 0}^n ( \boldsymbol{\chi}_{j+1} - \frac{\boldsymbol{C}_j}{j+1})\boldsymbol{\xi}_2$. The $0$-th term is
$\sim \boldsymbol{C}_0 \boldsymbol{\xi}_3 + 
\boldsymbol{C}_0 \boldsymbol{\xi}_2 \log n = o(n)$ and the martingale differences are bounded. \hfill $\Box$

It may be noted from the above that, as a consequence of the Azuma inequality for martingales with
uniformly bounded increments, one again gets  $\boldsymbol{C}_{n} \boldsymbol{\xi}_3/(n+1) \rightarrow 0$ almost surely.
 
Next for the repeated eigenvalue $\lambda \neq 0$  we consider the $0$-th term and the size of the martingale differences 
for $\boldsymbol{C}_{n+1} \boldsymbol{\xi}_3$ where $\boldsymbol{\xi}_3$ is the Jordan vector, the case of the eigenvector
$\boldsymbol{\xi}_2$ having been worked out in Theorem \ref{prop1}

\begin{proposition}\label{xi3}
For $\lambda \in (-1,1) \backslash \{0\}$, $\boldsymbol{C}_{n+1} \boldsymbol{\xi}_3$ when expanded in terms of martingale differences
has a $0$-th term which is $o(n)$ and its sum of squares of bounds on martingale differences is bounded by constant times 
$(1 + \log n)^2 D_n (\lambda)$.
\end{proposition}
 
{\bf Proof:} We now have 
\begin{eqnarray}\label{eqn4}
\boldsymbol{C}_{n+1} \boldsymbol{\xi}_3 &=& \boldsymbol{C}_n \boldsymbol{\xi}_3 + \boldsymbol{\chi}_{n+1} 
(\boldsymbol{\xi}_2 + \lambda \boldsymbol{\xi}_3) \nonumber \\
&=& (1 + \frac{\lambda}{n+1}) \boldsymbol{C}_n \boldsymbol{\xi}_3 + (\boldsymbol{\chi}_{n+1} - \frac{\boldsymbol{C}_n}{n+1})(\boldsymbol{\xi}_2 
+ \lambda \boldsymbol{\xi}_3) + 
\frac{\boldsymbol{C}_n}{n+1} \boldsymbol{\xi}_2 \nonumber \\
&=& \Big{(} \prod_0^n (1 + \frac{\lambda}{j+1}) \Big{)} \boldsymbol{C}_0 \boldsymbol{\xi}_3 \nonumber \\
&+& \sum_{j = 0}^n \Big{\{} \prod_{k = j+1}^n (1 + \frac{\lambda}{k+1}) \Big{\}} (\boldsymbol{\chi}_{j+1} - \frac{\boldsymbol{C}_j}{j+1})(\boldsymbol{\xi}_2 
+ \lambda \boldsymbol{\xi}_3) \nonumber \\
&+& \sum_{j = 0}^n \Big{\{} \prod_{k = j+1}^n (1 + \frac{\lambda}{k+1}) \Big{\}} 
 \frac{\boldsymbol{C}_j}{j+1} \boldsymbol{\xi}_2. 
\end{eqnarray}
The first two terms of the last line are familiar from the proof of Theorem \ref{prop1} (equation (\ref{martingale}) with appropriate changes), and the new last term 
\begin{eqnarray}\label{eqn5}
&& \sum_{j = 0}^n \Big{\{} \prod_{k = j+1}^n (1 + \frac{\lambda}{k+1}) \Big{\}} 
 \frac{\boldsymbol{C}_j}{j+1} \boldsymbol{\xi}_2.   \nonumber \\
&=& \sum_{j = 0}^n \Big{\{} \prod_{k = j+1}^n (1 + \frac{\lambda}{k+1}) \Big{\}} \nonumber \times \\
&&  \frac{1}{j+1} \Big{[} \Big{(} \prod_0^{j-1} (1 + \frac{\lambda}{k+1}) \Big{)} \boldsymbol{C}_0.\boldsymbol{\xi}_2 \nonumber \\
&+& \sum_{i = 0}^{j-1} \Big{\{} \prod_{l = i+1}^{j-1} (1 + \frac{\lambda}{l+1}) \Big{\}} 
\lambda (\boldsymbol{\chi}_{i+1} - \frac{\boldsymbol{C}_i}{i+1})\boldsymbol{\xi}_2  \Big{]},
\end{eqnarray}
using equation (\ref{martingale}) and the notational convention in the line following it.
Concentrating on (\ref{eqn5}) we first look at the new terms multiplying the martingale differences $ (\boldsymbol{\chi}_{i+1} - 
\frac{\boldsymbol{C}_i}{i+1})\boldsymbol{\xi}_2$  which are (after interchanging the order of summation)
\begin{equation}\label{eqn6}
\sum_{j = i + 1}^n \Big{\{} \prod_{k = j+1}^n (1 + \frac{\lambda}{k+1}) \Big{\}} 
\times \frac{1}{j+1} \Big{\{} \prod_{l = i+1}^{j-1} (1 + \frac{\lambda}{l+1}) \Big{\}}, 
\end{equation}
and for $i = 1, \ldots , n$ they are bounded above by constant times \\
$\sum_{j=i + 1}^n \frac{n^\lambda}{j^{\lambda +1}}  \frac{(j-1)^\lambda}{i^\lambda} \leq  \frac{n^\lambda}{i^\lambda}(1 + \log n) 
\max_{k \in \{2, \ldots, n\}} (1- \frac{1}{k})^\lambda$,
using $\sum_{j=i+1}^n \frac{1}{j} \leq \sum_{j = 1}^n \frac{1}{j} \leq 1 + \int_1^n \frac{1}{x} dx$, since $g(x) = \frac{1}{x}$ is decreasing.
Now, $\max_{k \in \{2, \ldots, n\}} (1- \frac{1}{k})^\lambda$ is bounded by 1 for $\lambda > 0$ and by $(1 - \frac{1}{2})^\lambda$ for $\lambda < 0$.
The term for $i = 0$ in (\ref{eqn6}) is also bounded by  constant times $n^\lambda (1 + \log n)$. 
Thus for the sum of squares of martingale difference bounds, instead of $D_n(\lambda)$ as in the previous application of Azuma-Hoeffding bound, we 
can now use the bound constant times $(1 + \log n)^2 D_n(\lambda)$. Next, the $0$-th term involving $\boldsymbol{C}_0$ in (\ref{eqn5}) is bounded by constant times 
$n^{\lambda} + n^{\lambda} + \sum_{j = 2}^n \frac{n^\lambda}{j^\lambda} \frac{1}{j} . (j-1)^\lambda \leq const. n^\lambda (1 + \log n),$
where for $\lambda < 0$ we use the upper bound $(1-\frac{1}{2})^\lambda$ for $(1 - \frac{1}{j})^\lambda$ (see the Appendix for the terms).

Finally, to get the required $0$-th term for $\boldsymbol{C}_{n+1} \boldsymbol{\xi}_3$ , we collect $0$-th terms of (\ref{eqn4}), and use the estimates from
Theorem \ref{prop1} and the previous paragraph. We also collect the martingale differences and 
use $(a + b)^2 \leq 2(a^2 + b^2)$ to get bounds on the sums of squares of the martingale difference bounds. \hfill $\Box$

The statement and proof of the Azuma-Hoeffding bound for $\alpha_2 \boldsymbol{C}_{n+1} . \boldsymbol{\xi_2} + \alpha_3 \boldsymbol{C}_{n+1} . \boldsymbol{\xi_3}$
where $\boldsymbol{\xi}_2, \boldsymbol{\xi}_3$ are the eigenvector and Jordan vector respectively to the repeated eigenvalue $\lambda
(\neq 0)$ is now similar to the statement of Theorem \ref{prop2} with the appropriate changes brought in by 
Proposition \ref{xi3}. We omit the detailed formulas, which can be recovered as needed.

{\bf Remark:} In cases with more colors, other Jordan vectors corresponding to the same real eigenvalue 
involve multiple iterated sums and are similar though more involved, and we have not pursued them in this short article. Also in
the case of complex eigenvalues a possible approach is to take the real part of the right hand side of equation (\ref{martingale}),
we refer to Basak and Dasgupta \cite{bd}, Dasgupta and Maulik 
\cite{dasguptamaulik}, Janson \cite{janson} etc. and the references therein. 

\section{Concluding remarks}

In strong or weak limit theorems for urn models it is usual to derive recursive equations for 
appropriately scaled $\boldsymbol{C}_{n+1} \boldsymbol{\xi}$, for example in law of large 
numbers the scaling is $(n+1)$, in central limit theorem the scaling is $\sqrt{n}$ or 
$\sqrt{n \log n}$ etc. depending on the eigenvalue. For large deviations 
we have used recursive equations for $\boldsymbol{C}_{n+1} \boldsymbol{\xi}$ itself
along with the the Azuma inequality to derive Azuma-Hoeffding bounds, which are
exponentially decreasing depending on the eigenvalues, for a class
of urn models through linear combinations corresponding to the eigenvectors and Jordan vectors
of the replacement matrix.

\section{Appendix} The coefficient of $\boldsymbol{C}_0 \boldsymbol{\xi}_2$ in equation (\ref{eqn5}) by direct calculation 
is 
\begin{eqnarray*}
&& \{ \prod_{k=1}^n (1 + \frac{\lambda}{k+1}) \} \times  \frac{1}{0 + 1} \times 1 + \{ \prod_{k = 2}^n(1 + \frac{\lambda}{k+1}) \} \times \frac{1}{1 + 1} \times (1 + \frac{\lambda}{1})\\
&+& \sum_{j = 2}^n \{ \prod_{k =j+1}^n (1 + \frac{\lambda}{k+1}) \} \times \frac{1}{j + 1} \times \{ \prod_{l = 0}^{j - 1} (1 + \frac{\lambda}{l+1}) \}.  
\end{eqnarray*}

\end{document}